\documentclass[10pt,a4paper]{amsart}

\usepackage{amssymb}
\usepackage{amsthm}
\usepackage{amsfonts}
\usepackage{microtype}
\usepackage{mathrsfs}
\usepackage{graphicx}
\usepackage{color}
\usepackage{latexsym}
\usepackage{rotating}
\usepackage{xspace}
\usepackage[all]{xy}
\usepackage{longtable}
\usepackage{leftidx}
\usepackage{mathtools}
\usepackage{titlesec}
\usepackage{tikz}
\usetikzlibrary{calc}
\usetikzlibrary{arrows}
\usetikzlibrary{decorations.markings}
\usepackage{geometry}

\newtheorem{theorem}{Theorem}[section]
\numberwithin{equation}{theorem}
\newtheorem{lemma}[theorem]{Lemma}
\newtheorem{proposition}[theorem]{Proposition}
\newtheorem{corollary}[theorem]{Corollary}

\theoremstyle{definition}
\newtheorem{definition}[theorem]{Definition}
\newtheorem{example}[theorem]{Example}
\newtheorem{notation}[theorem]{Notation}
\newtheorem{observation}[theorem]{Observation}

\theoremstyle{conjecture}

\newcommand{\Ass}{\operatorname{Ass}}

\newcommand{\grade}{\operatorname{grade}}

\newcommand{\Spec}{\operatorname{Spec}}

\newcommand{\ara}{\operatorname{ara}}

\newcommand{\cd}{\operatorname{cd}}

\newcommand{\id}{\operatorname{id}}

\newcommand{\pd}{\operatorname{pd}}

\newcommand{\V}{\operatorname{V}}
\newcommand{\RH}{\operatorname{H}}

\newcommand{\Ext}{\operatorname{Ext}}
\newcommand{\Supp}{\operatorname{Supp}}

\newcommand{\Tor}{\operatorname{Tor}}
\newcommand{\Hom}{\operatorname{Hom}}

\newcommand{\Ann}{\operatorname{Ann}}
\newcommand{\Rad}{\operatorname{Rad}}

\newcommand{\depth}{\operatorname{depth}}

\newcommand{\Max}{\operatorname{Max}}

\newcommand{\lo}{\longrightarrow}
\newcommand{\fm}{\frak{m}}
\newcommand{\fp}{\frak{p}}

\newcommand{\fa}{\frak{a}}
\newcommand{\fb}{\frak{b}}

\newcommand{\suchthat}{\;\ifnum\currentgrouptype=16 \middle\fi|\;}

\newenvironment{prf}[1][Proof]{\begin{proof}[\bf #1]}{\end{proof}}

\makeatletter
\newcommand{\holim@}[2]{%
\vtop{\m@th\ialign{##\cr
\hfil$#1\operator@font holim$\hfil\cr
\noalign{\nointerlineskip\kern1.5\ex@}#2\cr
\noalign{\nointerlineskip\kern-\ex@}\cr}}%
}
\newcommand{\holim}{%
\mathop{\mathpalette\holim@{\rightarrowfill@\textstyle}}\nmlimits@
}
\makeatother

\makeatletter
\def\@secnumfont{\bfseries}
\def\section{\@startsection{section}{1}%
\z@{.7\linespacing\@plus\linespacing}{.5\linespacing}%
{\normalfont\Large\bfseries\filcenter}}
\def\subsection{\@startsection{subsection}{2}%
\z@{.5\linespacing\@plus.7\linespacing}{-.5em}%
{\normalfont\large\bfseries}}
\makeatother
\begin{document}

\author[P. Pourghobadian, K. Divaani-Aazar and A. Rahimi]
{Parisa Pourghobadian, Kamran Divaani-Aazar and Ahad Rahimi}

\title[Relative homological ...]
{Relative homological rings and modules}

\address{P. Pourghobadian, Department of Mathematics, Faculty of Mathematical Sciences, Alzahra
University, Tehran, Iran.}
\email{paparpourghobadian@gmail.com}

\address{K. Divaani-Aazar, Department of Mathematics, Faculty of Mathematical Sciences, Alzahra
University, Tehran, Iran.}
\email{kdivaani@ipm.ir}

\address{A. Rahimi, Department of Mathematics, Razi University, Kermanshah, Iran.}
\email{ahad.rahimi@razi.ac.ir}

\subjclass[2020]{13C14; 13C05; 13D45.}

\keywords {Cohen-Macaulay module; cohomological dimension; complete intersection ring; Foxby equivalence; Gorenstein ring; local
cohomology; regular ring; relative Cohen-Macaulay module; relative dualizing module; relative system of parameters; semidualizing
module.}

\begin{abstract} The study of rings and modules with homological criteria is fundamental to commutative algebra. Consider a
commutative Noetherian ring $R$ with identity (which need not be local) and a proper ideal $\fa$ of $R$. In this paper, we
develop a relative analogue of the theory of homological rings and modules. Specifically, we introduce the notions of
$\fa$-relative regular, $\fa$-relative complete intersection, and $\fa$-relative Gorenstein rings and modules. By demonstrating
some interactions between these types of rings and modules, we extend some classical results.
\end{abstract}

\maketitle

\tableofcontents

\section{Introduction}

Throughout this article, the word ``ring" stands for a commutative Noetherian ring with identity. Let $\fa$ be a proper ideal
of a ring $R$. In the following, we will refer to a ring (resp. module) having a homological criterion as a ``homological ring"
(resp. ``homological module"). The theory of homological rings and modules traces back to 1954, when Auslander, Buchsbaum, and
Serre established a celebrated homological criterion for regular local rings. Since the 1960s, the study of homological
conjectures has been a significant research area in commutative algebra. Recently, Yves Andr\'{e} \cite{An} made a significant
breakthrough in these conjectures by using the theory of perfectoid spaces. Homological rings and modules are the focus of most
of these conjectures. Cohen-Macaulay modules are the most significant homological modules, and there are various generalizations
of them in the literature, including the notion of relative Cohen-Macaulay modules.

The theory of relative Cohen-Macaulay modules was first introduced by Hellus and Schenzel \cite{HeSc} and Rahro Zargar and Zakeri
\cite{RZ2}. A finitely generated $R$-module $M$ is said to be $\fa$-{\it relative Cohen-Macaulay} if $\RH_{\fa}^{i}\left(M\right)=0$
for all $i\neq \cd(\fa,M)$, where $\cd(\fa,M)$ denotes the {\it cohomological dimension} of $M$ with respect to $\fa$; that is, the
largest integer $i$ for which $\RH_{\fa}^{i}\left(M\right)\neq 0$. The study of relative Cohen-Macaulay modules has been pursued by
several authors; see, for instance, \cite{HeSt1, Sc1, Sc2, R, JR, Ra2, RZ1, CH, Ra1, DGTZ1, DGTZ2}.

Several subclasses of Cohen-Macaulay rings that have been a subject of research for many years are also homological. These are regular,
complete intersection, and Gorenstein rings, which satisfy the following implications: $\text{regular ring} \Rightarrow  \text{complete intersection ring} \Rightarrow \text{Gorenstein ring}.$

The aim of this paper is to establish a relative theory of regular, complete intersection, and Gorenstein rings and modules. We uncover
various interactions among these types of rings and modules, which expand some of the existing outcomes in the classical theory. We
demonstrate that certain results for homological modules do not hold for relative homological modules, as evidenced by several
counterexamples. The paper is organized as follows:

Section 2 provides several additional findings on $\fa$-relative Cohen-Macaulay modules, as shown in Propositions \ref{2.9d}, \ref{2.11f},
and \ref{2.12g}. Then, we introduce the concepts of $\fa$-relative regular modules, $\fa$-relative complete intersection rings, and
$\fa$-relative Gorenstein modules and compare them in Theorem \ref{2.19}.

Section 3 focuses on $\fa$-relative Gorenstein modules. Firstly, we introduce the notion of $\fa$-relative injective dimension of
$R$-modules to characterize $\fa$-relative Gorenstein modules. We establish that a finitely generated $R$-module $M$ with $M\neq
\fa M$ is $\fa$-relative Gorenstein if and only if $M$ is $\fa$-relative maximal Cohen-Macaulay and the $\fa$-relative injective dimension
of ${\RH}_{\fa}^{\cd(\fa,R)}(M)$ is zero; see Theorem \ref{3.2}. In Theorem \ref{3.5}, we present a technique for constructing
$\fa$-relative Gorenstein rings. We also provide some counterexamples to illustrate that the analogues of certain results for Gorenstein
rings do not hold for $\fa$-relative Gorenstein modules.

In Section 4, we investigate $\fa$-relative regular modules. We begin by assuming that $\fa$ is contained in the Jacobson radical of $R$,
and $M$ is a non-zero finitely generated $R$-module. If $M$ is $\fa$-relative regular, we establish in Theorem \ref{4.1a} that $$\pd_R(M/\fa M)=\pd_RM+\cd\left(\fa,R\right).$$ Then, in Theorem \ref{4.4d}, we demonstrate that $M$ is $\fa$-relative regular if and only if $$\grade(\fa,M)=\grade(\fa,R)=\mu(\fa),$$ where $\mu(\fa)$ denotes the minimum number of generators of $\fa$. Lastly, we establish the
invariance of the category of $\fa$-relative regular modules under certain equivalences and dualities of categories in Propositions
\ref{4.8h} and \ref{4.9i}.

\section{Relative Cohen-Macaulay modules}

Proposition \ref{2.9d} states that for any two nonzero finitely generated $R$-modules $M$ and $N$, if the ideal $\fa$ is contained
in the Jacobson radical of $R$ and $N$ is $\fa$-relative Cohen-Macaulay with $\cd(\fa,N)=\ara(\fa,N)$, then we have $$\grade(\fa,M)\leq \grade(\text{Ann}_RN,M)+\cd(\fa,N).$$
Proposition \ref{2.11f} provides a criterion for $\fa$-relative Cohen-Macaulay $R$-modules in terms of associated primes, assuming
$\fa$ is contained in the Jacobson radical of $R$. Next, in Proposition \ref{2.12g}, we show that if $R$ admits a faithful $\fa$-relative
Cohen-Macaulay module of finite projective dimension, then $R$ is $\fa$-relative Cohen-Macaulay. Lastly, after defining the notions
of $\fa$-relative regular modules, $\fa$-relative complete intersection rings and $\fa$-relative Gorenstein modules, we compare them
in Theorem \ref{2.19}.

For each non-negative integer $i$, the $i$-th local cohomology module of an $R$-module $M$ is defined as follows:
$$\RH_{\fa}^{i}\left(M\right)=\varinjlim \limits_{n\in \mathbb{N}} \text{Ext}_R^i\left(R/\fa^n,M\right).$$

To begin, we recall the definition of relative Cohen-Macaulay modules.

\begin{definition}\label{2.1} Let $\fa$ be a proper ideal of $R$ and $M$ a finitely generated $R$-module. Then $M$ is said to be
$\fa$-{\it relative Cohen-Macaulay} if either $M=\fa M$ or $M\neq \fa M$ and $\grade(\fa,M)=\cd(\fa,M)$.
\end{definition}

To prove Proposition \ref{2.9d}, we need the following preparatory results.

\begin{lemma}\label{2.2} Let $\fa$ be an ideal of $R$, $M$ an $\fa$-relative Cohen-Macaulay $R$-module with $M\neq \fa M$ and set
$c=\cd(\fa,M)$. Then
\begin{enumerate}
\item[(i)]  $\Supp_R({\RH}_{\fa}^c(M))=\Supp_R(M/\fa M)$.
\item[(ii)] $M_{\fp}$ is an $\fa R_{\fp}$-relative Cohen-Macaulay $R_{\fp}$-module and $\cd(\fa R_{\fp},M_{\fp})=c$ for every
$\fp\in \Supp_R(M/\fa M)$.
\item[(iii)] $M_{\fp}$ is a Cohen-Macaulay $R_{\fp}$-module for every minimal element $\fp$ of $\Supp_R(M/\fa M)$.
\end{enumerate}
\end{lemma}

\begin{prf} (i) It is obvious, since for every finitely generated $R$-module $L$, it is known and straightforward to verify
that $$\underset{i\in \mathbb{N}_0}\bigcup\Supp_R({\RH}_{\fa}^i(L))=\Supp_R(L/\fa L).$$

(ii) The flat base change theorem \cite[Theorem 4.3.2]{BS} implies an $R_{\fp}$-isomorphism ${\RH}_{\fa R_{\fp}}^j(M_{\fp})\cong {\RH}_{\fa}^j(M)_{\fp}$ for every prime ideal $\fp$ of $R$ and all $j\geq 0$.  Hence, as $M$ is $\fa$-relative Cohen-Macaulay,
it follows that ${\RH}_{\fa R_{\fp}}^j(M_{\fp})=0$ for every prime ideal $\fp$ of $R$ and all $j\neq c$. Let $\fp\in \Supp_R(M/\fa M)$.
Then (i) yields that ${\RH}_{\fa R_{\fp}}^c(M_{\fp})\neq 0$. Therefore, $$\grade(\fa R_{\fp},M_{\fp})=c=\cd(\fa R_{\fp},M_{\fp}),$$
and so $M_{\fp}$  is an $\fa R_{\fp}$-relative Cohen-Macaulay $R_{\fp}$-module.

(iii) Let $\fp$ be a minimal element of $\Supp_R(M/\fa M)$. Since $\fp$ is minimal over $\fa+\Ann_RM$, it follows that $$\Rad(\fa R_{\fp}+\Ann_{R_{\fp}}M_{\fp})=\fp R_{\fp},$$ and so Grothendieck's non-vanishing theorem \cite[Theorem 6.1.4]{BS} asserts that $$\dim_{R_{\fp}}M_{\fp}=\cd(\fp R_{\fp},M_{\fp})=\cd(\fa R_{\fp}+\Ann_{R_{\fp}}M_{\fp},M_{\fp})=\cd(\fa R_{\fp},M_{\fp}).$$
Hence,
$$
\begin{array}{ll}
\grade(\fa,M) & \leq \grade(\fa R_{\fp},M_{\fp})\\
&\leq \depth_{R_{\fp}}M_{\fp}\\
&\leq \dim_{R_{\fp}}M_{\fp}\\
& =\cd(\fa R_{\fp},M_{\fp})\\
& \leq \cd(\fa,M)\\
&= \grade(\fa,M).
\end{array}
$$
Thus $\depth_{R_{\fp}}M_{\fp}=\dim_{R_{\fp}}M_{\fp}$, and so $M_{\fp}$ is a Cohen-Macaulay $R_{\fp}$-module.
\end{prf}

\begin{lemma}\label{2.3}  Let $\fa$ be a proper ideal of $R$ and $M$ a nonzero finitely generated $R$-module. Consider the following
conditions:
\begin{enumerate}
\item[(i)] $M$ is $\fa$-torsion.
\item[(ii)]  $\cd(\fa,M)=0$.
\item[(iii)]  $\Ass_RM\cap \V(\fa)\neq \emptyset$.
\item[(iv)]  $\Ass_RM\subseteq \V(\fa)$.
\end{enumerate}
Then  (i) $\Rightarrow$ (ii), (i) $\Leftrightarrow$ (iv) and (iv) $\Rightarrow$ (iii) hold. If $\fa$ is contained in the Jacobson
radical of $R$, then (ii) $\Rightarrow$ (i) holds. Also, (iii) $\Rightarrow$ (ii) holds provided that $M$ is $\fa$-relative
Cohen-Macaulay. So, these four conditions are equivalent when $\fa$ is contained in the Jacobson radical of $R$ and $M$ is
$\fa$-relative Cohen-Macaulay.
\end{lemma}

\begin{prf} The implications (i) $\Rightarrow$ (ii), (i) $\Leftrightarrow$ (iv), and (iv) $\Rightarrow$ (iii) are known and clear.

Now, suppose that $\fa$ is contained in the Jacobson radical of $R$ and (ii) holds. Set $\widetilde{M}=M/\Gamma_{\fa}(M)$. Then,
${\RH}_{\fa}^0(\widetilde{M})=0$ and ${\RH}_{\fa}^i(\widetilde{M})\cong {\RH}_{\fa}^i(M)=0$ for all $i>0$. As we mentioned in the
proof of Lemma \ref{2.2}(i), for every finitely generated $R$-module $L$, one has $$\underset{i\in \mathbb{N}_0}\bigcup
\Supp_R({\RH}_{\fa}^i(L))=\Supp_R(L/\fa L).$$ Hence, $\Supp_R(\widetilde{M}/\fa \widetilde{M})=\emptyset$, and so $\widetilde{M}=
\fa \widetilde{M}$. Thus $M=\Gamma_{\fa}(M)$ by Nakayama's lemma, and so (i) holds.

Finally, assume that $M$ is $\fa$-relative Cohen-Macaulay and (iii) holds. Let $\fp\in \Ass_RM\cap \V(\fa)$. Then $$\grade(\fa R_{\fp},
M_{\fp})\leq \depth_{R_{\fp}}M_{\fp}=0.$$ Now, by Lemma \ref{2.2}(ii), the $R_{\fp}$-module $M_{\fp}$ is $\fa R_{\fp}$-relative
Cohen-Macaulay and $$\cd(\fa,M)=\cd(\fa R_{\fp},M_{\fp})=\grade(\fa R_{\fp},M_{\fp})=0,$$ and so (ii) holds.
\end{prf}

\begin{definition}\label{2.4} Let $M$ be a finitely generated $R$-module and $\fa$ an ideal of $R$ with $M\neq \fa M$.
\begin{enumerate}
\item[(i)] Let $c=\cd\left(\fa,M\right)$. A sequence $x_{1}, x_2, \ldots, x_{c}\in \fa$ is called $\fa$-{\it relative system of
parameters}, $\fa$-s.o.p, of $M$ if $$\Rad\left(\langle x_{1}, x_2, \ldots, x_{c}\rangle+\Ann_{R}M\right)=\Rad\left(\fa+\Ann_{R}M
\right).$$
\item[(ii)] {\it Arithmetic rank} of $\fa$ with respect to $M$, $\ara\left(\fa,M \right)$, is defined as the infimum of the integers
$n\in \mathbb{N}_0$ such that there exist $x_1, x_2, \ldots, x_n\in R$ satisfying $$\Rad\left(\langle x_{1}, x_2, \ldots, x_{n}\rangle +\Ann_{R}M\right)=\Rad\left(\fa+\Ann_{R}M\right).$$
\end{enumerate}
\end{definition}

Clearly, if $x_{1},x_2,\ldots, x_{c}\in R$ is an $\fa$-s.o.p of $M$, then for all $t_{1},\ldots,t_{c}\in \mathbb{N}$, every permutation
of $x_{1}^{t_{1}},\ldots, x_{c}^{t_{c}}$ is also an $\fa$-s.o.p of $M$. One may easily check that $\cd\left(\fa,M\right)\leq \ara
\left(\fa,M \right)$. Obviously, $\ara\left(\fa,R\right)=\ara\left(\fa\right)$.

Let $M$ be a $d$-dimensional finitely generated module over a local $(R,\fm)$ and $x_1,\dots, x_d\in \fm$. It is evident that
$x_1,\dots, x_d$ is a system of parameters of $M$ if and only if $x_1,\dots, x_d$ is an $\fm$-relative system of parameters of
$M$. In spite of the fact that $M$ always admits an $\fm$-relative system of parameters, this is not true for a general ideal
$\fa$. In this regard, we have the following:

\begin{lemma}\label{2.5} (See \cite[Lemma 2.2]{DGTZ2}.) Let $M$ be a finitely generated $R$-module and $\fa$ an ideal of $R$ with
$M\neq \fa M$. Then $\fa$ contains an $\fa$-s.o.p of $M$ if and only if $\cd\left(\fa,M\right)=\ara\left(\fa,M \right)$.
\end{lemma}

\begin{lemma}\label{2.6a} (See \cite[Lemma 2.4 and Theorem 2.7]{DGTZ2}.) Let $\fa$ be an ideal of $R$, $M$ a finitely generated
$R$-module with $M\neq \fa M$ and $c=\cd\left(\fa,M\right)$. Assume that $\cd\left(\fa,M\right)=\ara\left(\fa,M \right)$ and $x_1,
\ldots, x_{c}\in \fa$. Consider the following conditions:
\begin{enumerate}
\item[(i)] $x_1,\ldots, x_{c}$ is an $\fa$-s.o.p of $M$.
\item[(ii)] $\cd\left(\fa,M/\langle x_{1},x_2, \ldots, x_{i}\rangle M\right)=c-i$ for every $i=1, 2,\ldots, c$.
\end{enumerate}
Then (i) implies (ii). Additionally, if $\fa$ is contained in the Jacobson radical of $R$, then (i) and (ii) are equivalent.
\end{lemma}

\begin{lemma}\label{2.7b} (See \cite[Theorem 3.3]{DGTZ2}.) Let $M$ be a finitely generated $R$-module and $\fa$ an ideal of $R$ with $\cd\left(\fa,M\right)=\ara\left(\fa,M \right)$. Consider the following conditions:
\begin{enumerate}
\item[(i)] $M$ is $\fa$-relative Cohen-Macaulay.
\item[(ii)] Every $\fa$-s.o.p of $M$ is an $M$-regular sequence.
\item[(iii)] There exists an $\fa$-s.o.p of $M$ which is an $M$-regular sequence.
\end{enumerate}
Then (i) and (iii) are equivalent. Furthermore, if $\fa$ is contained in the Jacobson radical of $R$, then all three conditions
are equivalent.
\end{lemma}

In light of Lemma \ref{2.7b}, one may wonder whether every $\fa$-relative Cohen-Macaulay module has an $\fa$-s.o.p. The following
example shows that this is not the case.

\begin{example}\label{2.8c} Let $\Bbbk$ be a field and $S=\Bbbk[[x,y,z,w]]$. Consider the elements $f=xw-yz$, $g=y^{3}-x^{2}z$,
and $h=z^{3}-y^2w$ of $S$. Let $R=S/\langle f \rangle$, and $\mathfrak{a}=\langle f,g,h \rangle/\langle f \rangle$. Then $R$ is
a local complete intersection ring of dimension $3$, $\cd(\mathfrak{a},R)=1$, and $\ara(\mathfrak{a})\geq 2$; see
\cite[Remark 2.1(ii)]{HeSt2}. As $\cd(\mathfrak{a},R)\neq \ara(\mathfrak{a})$, by Lemma \ref{2.5}, $R$ possesses no $\fa$-s.o.p.
Since ${\RH}_{\fa}^0(R)=0$, it follows that $R$ is $\fa$-relative Cohen-Macaulay.
\end{example}

The next result concerns the vanishing of certain Ext modules.

\begin{proposition}\label{2.9d} Let $\fa$ be an ideal of $R$ contained in its Jacobson radical and $M$ and $N$ two nonzero finitely
generated $R$-modules.  Assume that $N$ is $\fa$-relative Cohen-Macaulay and $\cd(\fa,N)=\ara(\fa,N)$. Then $\Ext_R^i(N,M)=0$ for all $i<\grade(\fa,M)-\cd(\fa,N)$.
\end{proposition}

\begin{prf} We proceed by induction on $c=\cd(\fa,N)$. If $c=0$, then by Lemma \ref{2.3}, we conclude that $\Supp_RN\subseteq \V(\fa)$,
and so the claim in this case follows by Rees' theorem; see e.g. \cite[Theorem 16.6]{M}.

Next suppose that $c>0$, and let $x_1, \dots, x_c\in \fa$ be an $\fa$-s.o.p of $N$. Set $\overline{N}=N/x_1N$. Then Lemma \ref{2.6a}
implies that $\cd(\fa,\overline{N})=c-1$, and $x_1, \dots, x_c$ is an $N$-regular sequence by Lemma \ref{2.7b}. As $\grade(\fa,\overline{N})=
\grade(\fa,N)-1$, it turns out that $\overline{N}$ is an $\fa$-relative Cohen-Macaulay $R$-module. One can easily check that $x_2, \dots,
x_c\in \fa$ is an $\fa$-s.o.p of $\overline{N}$, and so $\cd(\fa,\overline{N})=\ara(\fa,\overline{N})$ by Lemma \ref{2.5}.
From the short exact sequence $$0\lo N\overset{x_1}\lo N\lo \overline{N}\lo 0,$$ one obtains the exact sequence $$\cdots \lo \Ext_R^{i}(\overline{N},M)\lo \Ext_R^{i}(N,M)\overset{x_1}\lo \Ext_R^{i}(N,M)\lo \Ext_R^{i+1}(\overline{N},M)\lo \cdots .$$
Let $i<\grade(\fa,M)-c$ be an integer. Then $i+1<\grade(\fa,M)-(c-1)$, and so $\Ext_R^{i+1}(\overline{N},M)=0$ by the induction
hypothesis. Thus, the map $\Ext_R^{i}(N,M)\overset{x_1}\lo \Ext_R^{i}(N,M)$ is surjective, and so $\Ext_R^{i}(N,M)=0$ by Nakayama's lemma.
\end{prf}

The following corollary provides a lower bound for the vanishing of generalized local cohomology modules. Recall that, for an ideal $\fa$
of $R$ and two $R$-modules $N$ and $M$, the $i$-th generalized local cohomology module of $N$ and $M$ with respect to $\fa$ is defined by $$\RH^i_{\fa}(N,M)=\underset{n\in \mathbb{N}}{\varinjlim} \ \text{Ext}^{i}_{R}(N/\fa^{n}N,M);$$ see \cite{H}.

\begin{corollary}\label{2.10e} Let $\fa$ be an ideal of $R$ contained in its Jacobson radical and $M$ and $N$ two nonzero finitely generated
$R$-modules.  Assume that $N$ is $\fa$-relative Cohen-Macaulay and $\cd(\fa,N)=\ara(\fa,N)$. Then $\RH_{\fb}^n(N,M)=0$ for every ideal
$\fb$ of $R$ and all $n<\grade(\fa,M)-\cd(\fa,N)$.
\end{corollary}

\begin{prf} By the proof of \cite[Theorem 2.5]{DH}, there is the following Grothendieck spectral sequence
$$\text{E}^{p,q}_2:=\RH^{p}_{\fb} (\Ext_R^q(N,M))\underset{p}\Longrightarrow \RH^{p+q}_{\fb}(N,M).$$ So, for each non-negative integer $n$,
there exists a chain $$0=\RH^{-1}\subseteq \RH^{0}\subseteq \cdots \subseteq \RH^{n}=\RH^{n}_{\fb}(N,M) \  \  (\ast) $$ of submodules of $\RH^{n}_{\fb}(N,M)$ such that $\RH^{p}/\RH^{p-1}\cong \text{E}^{p,n-p}_{\infty}$ for all $p=0,1,\dots, n$. Let $p$ and $n$ be two integers
such that $$0\leq p\leq n<\grade(\fa,M)-\cd(\fa,N).$$ Then $\Ext_R^{n-p}(N,M)=0$ by Proposition \ref{2.9d}. This implies that $\text{E}^{p,n-p}_{\infty}=0$, because $\text{E}^{p,n-p}_{\infty}$ is a subquotient of $\text{E}^{p,n-p}_2$. Therefore, from the chain $(*)$,
we can deduce that $\RH_{\fb}^n(N,M)=0$.
\end{prf}

Next, we present a new characterization of $\fa$-relative Cohen-Macaulay modules in the case that $\fa$ is contained in the Jacobson radical
of $R$.

\begin{proposition}\label{2.11f} Let $\fa$ be an ideal of $R$ contained in its Jacobson radical and $M$ a nonzero finitely generated
$R$-module. Then the following are equivalent:
\begin{enumerate}
\item[(i)] $M$ is $\fa$-relative Cohen-Macaulay;
\item[(ii)] $\cd(\fa,R/\fp)=\grade(\fa,M)$ for all $\fp\in \Ass_RM$.
\end{enumerate}
\end{proposition}

\begin{prf} (i) $\Rightarrow$ (ii) First, assume that $\cd(\fa,M)=0$. Then $\grade(\fa,M)=0$. On the other hand, \cite[Theorem 2.2]{DNT}
yields that $$0\leq \cd(\fa,R/\fp)\leq \cd(\fa,M)=0$$ for all $\fp\in \Ass_RM$.

Next, suppose that  $\cd(\fa,M)>0$. Then by Lemma \ref{2.3}, we deduce that $\Ass_RM=\Ass_RM\setminus \V(\fa)$. Therefore, by
\cite[Lemma 3.3]{DGTZ1}, we conclude that $$\cd(\fa,R/\fp)=\cd(\fa,M)=\grade(\fa,M)$$ for all $\fp\in \Ass_RM$.

(ii) $\Rightarrow$ (i)  Set $N=\underset{\fp\in \Ass_RM}\bigoplus R/\fp$. Then we can easily see that $\Supp_RN=
\Supp_RM$. Hence, by \cite[Theorem 2.2]{DNT}, we deduce that $$\cd(\fa,M)=\cd(\fa,N)=\max\{\cd(\fa,R/\fp)\mid \fp\in \Ass_RM \}=
\grade(\fa,M),$$ as desired.
\end{prf}

When a local ring $(R,\fm)$ admits a nonzero Cohen-Macaulay module of finite projective dimension, it follows that $R$ is Cohen-Macaulay
by the Peskine-Szpiro intersection theorem. We now present a partial relative analogue of this result.

\begin{proposition}\label{2.12g} Let $\fa$ be a proper ideal of $R$. Assume that $R$ admits a faithful $\fa$-relative Cohen-Macaulay
module $M$ of finite projective dimension. Then $R$ is $\fa$-relative Cohen-Macaulay.
\end{proposition}

\begin{prf} Set $c=\cd(\fa,R)$. For every finitely generated $R$-module $N$, \cite[Proposition 1.2.10(a)]{BH} implies that $$\grade(\fa,N)=\inf \{\depth N_{\fp} \mid \fp\in \V(\fa) \}.$$ Since $M$ is faithful, it follows that $\Supp_RM=\Spec R$, and so
$\cd(\fa,M)=c$ by \cite[Theorem 2.2]{DNT}. Let $\fp \in \V(\fa)$. As $\pd_RM<\infty$, it follows that $\pd_{R_{\fp}}M_{\fp}<\infty$, and so by
the Auslander-Buchsbaum formula, we get that $\depth M_{\fp}\leq \depth R_{\fp}.$ Therefore, $\grade(\fa,M)\leq \grade(\fa,R)$. Now, we have
$$\begin{array}{ll}
c&=\grade(\fa,M)\\
&\leq \grade(\fa,R)\\
&\leq \cd(\fa,R)\\
&=c,
\end{array}$$
and so $\grade(\fa,R)=\cd(\fa,R)$, as desired.
\end{prf}

\begin{definition}\label{2.13h} Let $\fa$ be a proper ideal of $R$.
\begin{enumerate}
\item[(i)] A finitely generated $R$-module $M$ is called $\fa$-{\it relative regular} if either $M=\fa M$ or $M\neq \fa M$ and
$\fa$ can be generated by an $M$-regular sequence of length $\cd(\fa,R)$ which is also an $R$-regular sequence.
\item[(ii)] Let $\mathcal{R}$ denote the $\fa$-adic completion of $R$. We say that $R$ is $\fa$-relative complete intersection
if there exist a ring $T$, a proper ideal $\fb$ of $T$, and elements $x_1, x_2, \dots, x_{\ell}\in \fb$ such that:

(1) $x_1, x_2, \dots, x_{\ell}$ is both a part of a $\fb$-s.o.p of $T$ and a $T$-regular sequence.

(2) $T$ is $\fb$-relative regular, $T/\langle x_1, x_2, \dots, x_{\ell} \rangle \cong \mathcal{R}$, and $\fb
\mathcal{R}=\fa \mathcal{R}$.

\item[(iii)] A finitely generated $R$-module $M$ is called $\fa$-{\it relative Gorenstein} if $\Ext_R^i(R/\fa,M)=0$ for all
$i\neq \cd(\fa,R)$.
\item[(iv)] A finitely generated $R$-module $M$ is called $\fa$-{\it relative maximal Cohen-Macaulay} if $\grade(\fa,M)=\cd(\fa,R)$.
\end{enumerate}
\end{definition}

Note that if there exists an $\fa$-relative regular $R$-module $M$ with $M\neq \fa M$, it implies that the ideal $\fa$ is a complete
intersection. In particular, the ring $R$ is $\fa$-relative regular if and only if $\fa$ is a complete intersection. Therefore, if
there exists an $\fa$-relative regular $R$-module $M$ with $M\neq \fa M$, it follows that the ring $R$ itself is $\fa$-relative regular.

In the following immediate observation, the connection between classical homological modules and relative homological modules is
revealed.

\begin{observation}\label{2.14i} Let $(R,\fm)$ be a local ring and $M$ a finitely generated $R$-module.
\begin{enumerate}
\item[(i)] $R$ is regular if and only if $R$ is $\fm$-relative regular.
\item[(ii)] $R$ is complete intersection if and only if $R$ is $\fm$-relative complete intersection.
\item[(iii)] $M$ is Gorenstein if and only if $M$ is $\fm$-relative Gorenstein.
\item[(iv)] $M$ is maximal Cohen-Macaulay if and only if $M$ is $\fm$-relative maximal Cohen-Macaulay.
\item[(v)] $M$ is Cohen-Macaulay if and only if $M$ is $\fm$-relative Cohen-Macaulay.
\end{enumerate}
\end{observation}

Theorem \ref{2.19} is the main result of this section. To prove it, we need the following four lemmas.

From Definition \ref{2.13h}, we have the following immediate result:

\begin{lemma}\label{2.15j} Let $\fa$ be a proper ideal of $R$ and $\mathcal{R}$ denote the $\fa$-adic completion of $R$. Then $R$ is an $\fa$-relative complete intersection if and only if $\mathcal{R}$ is an $\fa \mathcal{R}$-relative complete intersection.
\end{lemma}

\begin{lemma}\label{2.16k} Let $\fa$ be an ideal of $R$, $f:R\lo T$ a flat ring homomorphism and $M$ a finitely generated
$R$-module with $M\neq \fa M$. We have the following inequalities, with equality when $f$ is faithfully flat.
\begin{enumerate}
\item[(i)] $\grade(\fa,M)\leq \grade(\fa T,M\otimes_RT)$.
\item[(ii)] $\cd\left(\fa,M\right)\geq \cd\left(\fa T,M\otimes_RT\right)$.
\end{enumerate}
\end{lemma}

\begin{prf} (i)  For every non-negative integer $i$, there is a natural $T$-isomorphism $$\Ext_T^i(T/\fa T,M\otimes_RT)\cong \Ext_R^i(R/\fa,M)\otimes_RT.$$ In particular, $\Ext_T^i(T/\fa T,M\otimes_RT)\neq 0$ implies that $\Ext_R^i(R/\fa,M)\neq 0$,
and so $$\grade(\fa,M)\leq \grade(\fa T,M\otimes_RT).$$ If $f$  is faithfully flat, then $\Ext_T^i(T/\fa T,M\otimes_RT)\neq 0$
if and only if $\Ext_R^i(R/\fa,M)\neq 0$, and so $\grade(\fa,M)=\grade(\fa T,M\otimes_RT)$.

(ii) The flat base change theorem yields a natural $T$-isomorphism $${\RH}_{\fa T}^{i}\left(M\otimes_RT\right)\cong {\RH}_{\fa}^{i}\left(M\right)\otimes_RT$$ for all $i\geq 0$. In particular, ${\RH}_{\fa T}^{i}\left(M\otimes_RT\right)\neq 0$
implies that ${\RH}_{\fa}^{i}\left(M\right)\neq 0$, and so $\cd\left(\fa,M\right)\geq \cd\left(\fa T,M\otimes_RT\right)$. If $f$
is faithfully flat, then  ${\RH}_{\fa T}^{i}\left(M\otimes_RT\right)\neq 0$ if and only if ${\RH}_{\fa}^{i}\left(M\right)\neq 0$,
and so $\cd\left(\fa,M\right)=\cd\left(\fa T,M\otimes_RT\right)$.
\end{prf}

\begin{lemma}\label{2.17l} Let $\fa$ be an ideal of $R$ contained in its Jacobson radical and $\mathcal{R}$ denote the $\fa$-adic
completion of $R$. Let $M$ be a finitely generated $R$-module.
\begin{enumerate}
\item[(i)] If $M$ is $\fa$-relative regular, then $M\otimes_R\mathcal{R}$ is $\fa \mathcal{R}$-relative regular.
\item[(ii)] $M$ is $\fa$-relative Gorenstein if and only if $M\otimes_R\mathcal{R}$ is $\fa \mathcal{R}$-relative Gorenstein.
\end{enumerate}
\end{lemma}

\begin{prf} Let $c=\cd\left(\fa,R\right)$ and $\psi:R\lo \mathcal{R}$ denote the natural ring monomorphism. Since $\fa$ is contained
in the Jacobson radical of $R$, it follows that $\mathcal{R}$ is a faithfully flat $R$-algebra. Clearly, we may and do assume that
$M\neq \fa M$, and so $M\otimes_R\mathcal{R}\neq (\fa \mathcal{R}) (M\otimes_R\mathcal{R})$. By Lemma \ref{2.16k}(ii), one has
$\cd\left(\fa \mathcal{R},\mathcal{R}\right)=\cd\left(\fa,R\right)$.

(i)  Assume that $M$ is $\fa$-relative regular. Then, $\fa$ has generators $x_{1}, x_2, \dots, x_{c}$ which form both an $M$-regular
sequence and an $R$-regular sequence. This immediately yields that $\fa \mathcal{R}=\langle \psi(x_{1}), \psi(x_2), \ldots,
\psi(x_{c})\rangle_\mathcal{R}$. Clearly, $\psi(x_{1}), \psi(x_2), \ldots, \psi(x_{c})$ is both an $M\otimes_R\mathcal{R}$-regular
sequence and an $\mathcal{R}$-regular sequence. Thus, $M\otimes_R\mathcal{R}$ is $\fa \mathcal{R}$-relative regular.

(ii)  For every non-negative integer $j$, in view of the faithfulness of $\psi$ and the existence of the natural $\mathcal{R}$-isomorphism $$\Ext_R^j(R/\fa,M)\otimes_R\mathcal{R}\cong \Ext_\mathcal{R}^j(\mathcal{R}/\fa \mathcal{R},M\otimes_R\mathcal{R}),$$ it turns out that $\Ext_R^j(R/\fa,M)=0$ if and only if $\Ext_\mathcal{R}^j(\mathcal{R}/\fa \mathcal{R},M\otimes_R\mathcal{R})=0.$ This completes the proof
(ii).
\end{prf}

\begin{lemma}\label{2.18} Let $\fa$ be a proper ideal of $R$. If $R$ is $\fa$-relative regular, then $\cd\left(\fa,R\right)=\pd_R(R/\fa)$.
\end{lemma}

\begin{prf}  Let $c=\cd\left(\fa,R\right)$. Assume that $R$ is $\fa$-relative regular. Then, there is an $R$-regular sequence
$x_{1}, x_2, \dots, x_{c}$ that generates $\fa$. So, $$\pd_R(R/\fa)=\pd_R\left(R/\langle x_{1}, x_2, \dots, x_{c}\rangle\right)=c.$$
\end{prf}

It is now time to present the main result of this section.

\begin{theorem}\label{2.19} Let $\fa$ be a proper ideal of $R$.
\begin{enumerate}
\item[(i)] Suppose that $\fa$ is contained in the Jacobson radical of $R$. If $R$ is $\fa$-relative regular, then it is $\fa$-relative
complete intersection.
\item[(ii)] Suppose that $\fa$ is contained in the Jacobson radical of $R$. If $R$ is $\fa$-relative complete intersection, then it
is $\fa$-relative Gorenstein.
\item[(iii)]  Every $\fa$-relative Gorenstein module $M$ with $M\neq \fa M$ is $\fa$-relative maximal Cohen-Macaulay.
\item[(iv)] Every $\fa$-relative maximal Cohen-Macaulay module is $\fa$-relative Cohen-Macaulay.
\end{enumerate}
\end{theorem}

\begin{prf} (i) By Lemmas \ref{2.17l}(i) and \ref{2.15j}, we may and do assume that $R$ is $\fa$-adically complete. So, the claim is
immediate by the definition.

(ii) In view of Lemmas \ref{2.15j} and \ref{2.17l}(ii), we may and do assume that $R$ is $\fa$-adically complete. So, there
are a ring $T$, a proper ideal $\fb$ of $T$ and $x_1, x_2, \dots, x_{\ell}\in \fb$ which form both a part of a $\fb$-s.o.p
of $T$ and a $T$-regular sequence such that $T$ is $\fb$-relative regular, $T/\langle x_1, x_2, \dots, x_{\ell} \rangle\cong R$ and
$\fb R=\fa$.

Set $c=\cd(\fb,T)$. As $T$ is $\fb$-relative regular, the ideal $\fb$ can be generated by a $T$-regular sequence of length $c$.
Hence $c\leq \grade(\fb,T)\leq c$, and so $\grade(\fb,T)=c$. The fact that $T$ is $\fb$-relative regular, also implies that
$\pd_T(T/\fb)=c$ by Lemma \ref{2.18}. Thus, $\Ext_T^i(T/\fb,T)=0$ for all $i\neq c$. By Lemma \ref{2.6a}, we have $$\cd(\fa,R)
=\cd(\fb,T/\langle x_1, x_2, \dots, x_{\ell} \rangle)=\cd(\fb,T)-\ell=c-\ell. \  \
\   \  (\dag)$$ Clearly, $R/\fa\cong T/\fb$. Since $x_1, x_2, \dots, x_{\ell}$ is a $T$-regular sequence, by
\cite[\S 18, Lemma 2(i)]{M}, we have a $T$-isomorphism $$\Ext_R^n(R/\fa,R)\cong \Ext_T^{n+\ell}(T/\fb,T)\  \  \  \
(\ddag)$$ for all $n\geq 0$. Now, $(\dag)$ and $(\ddag)$ imply that $\Ext_R^{i}(R/\fa,R)=0$ for
all $i\neq \cd(\fa,R)$, and so $R$ is $\fa$-relative Gorenstein.

(iii) Let $M$ be an $\fa$-relative Gorenstein module with $M\neq \fa M$. Then $\Ext_R^{i}(R/\fa,M)\neq 0$ if and only
if $i=\cd(\fa,R)$. This means that $\grade(\fa,M)=\cd(\fa,R)$.

(iv) Let $M$ be an $\fa$-relative maximal Cohen-Macaulay $R$-module. Then $\grade(\fa,M)=\cd(\fa,R)$. Now, we have
$$\cd(\fa,R)=\grade(\fa,M)\leq \cd(\fa,M)\leq \cd(\fa,R),$$ and so $\grade(\fa,M)=\cd(\fa,M)$.
\end{prf}

We conclude this section with the following three examples. The first example provides a relative Cohen-Macaulay
$R$-module, which is not relative Gorenstein.

\begin{example}\label{2.20} Let $\Bbbk$ be a field and $G$ the following cyclic graph:
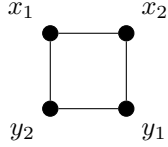
\begin{center}
\begin{tikzpicture}
\node[circle, draw, fill=black, inner sep=0pt, minimum size=2mm] (1) [label=south west:$y_2$] at (0,0) {};
\node[circle, draw, fill=black, inner sep=0pt, minimum size=2mm] (2) [label=north west:$x_1$] at (0,1) {};
\node[circle, draw, fill=black, inner sep=0pt, minimum size=2mm] (3) [label=north east:$x_2$] at (1,1) {};
\node[circle, draw, fill=black, inner sep=0pt, minimum size=2mm] (4) [label=south east:$y_1$] at (1,0) {};
\draw (1) -- (2) -- (3) -- (4) -- (1);
\end{tikzpicture}
\end{center}
Then the edge ideal of $G$ in the polynomial ring $S=\Bbbk[x_1,x_2,y_1,y_2]$ is $I(G)=\langle x_1x_2,x_2y_1,y_1y_2,y_2x_1 \rangle.$
Set $\fa=\langle y_1,y_2 \rangle$. It is routine to see that $I(G)=\langle x_1, y_1 \rangle \cap \langle x_2, y_2 \rangle$ is a minimal
primary decomposition of $I(G)$, and so $$\Ass_S(S/I(G))=\{\langle x_1, y_1 \rangle, \langle x_2, y_2 \rangle\}.$$ Since $\fa$ is not
contained in any member of $\Ass_S(S/I(G))$, it follows that $\grade(\fa,S/I(G))\geq 1$. On the other hand, one has $$\cd(\fa,S/I(G))=
\max \{\cd(\fa,S/\langle x_1, y_1 \rangle ),\cd(\fa,S/\langle x_2, y_2 \rangle)\}=1.$$ Hence, the $S$-module $S/I(G)$ is $\fa$-relative Cohen-Macaulay. But, $S/I(G)$ is not $\fa$-relative maximal Cohen-Macaulay, because $\cd(\fa,S)=2$. Consequently, $S/I(G)$ is not a
$\fa$-relative Gorenstein.
\end{example}

The following two examples show that there are plenty of relative Gorenstein and relative regular $R$-modules.

\begin{example}\label{2.21} Let $S=\Bbbk[x_1, \dots, x_n]$ be a polynomial ring over a field $\Bbbk$. Set $\fm=\langle x_1, \dots,
x_n\rangle$ and $\fa=\langle x_1^{p_1}, \dots, x_n^{p_n}, m_1, \dots, m_r\rangle,$ where $r, p_1,\dots, p_n\in\mathbb{N}$ and
$m_j$$'$s are monomials in $S$. By \cite[Proposition 28]{Bo}, we have $\Ext^i_S(S/\fa,S)=0$ for all $i\neq n$. Notice that
$\Rad(\fa)=\fm$. Hence, ${\RH}^i_\fa(S)={\RH}^i_{\fm}(S)$ for all $i\geq 0$, and so $\cd(\fa,S)=n$. Consequently, $S$ is
$\fa$-relative Gorenstein.
\end{example}

\begin{example}\label{2.22} Let $S=\Bbbk[x_1, \dots, x_n]$ be a polynomial ring over a field $\Bbbk$. Set $\fa=\langle
x_1^{p_1}, \dots, x_{\ell}^{p_{\ell}}\rangle,$ where $1\leq \ell\leq n$ and $p_1,\dots, p_{\ell}\in \mathbb{N}$. Then, obviously,
$S$ is $\fa$-relative regular.
\end{example}

\section{Relative Gorenstein modules}

In this section, we delve deeper into the study of relative Gorenstein modules. We begin by providing a characterization of relative
Gorenstein modules in Theorem \ref{3.2}. We then establish a class of relative Gorenstein rings in Theorem \ref{3.5}. However, we
also present some counterexamples to address some natural questions about relative Gorenstein modules.

It is well known that a $d$-dimensional local ring $(R,\fm)$ is Gorenstein if and only if it is Cohen-Macaulay and ${\RH}_{\fm}^d(R)\cong \text{E}_R(R/\fm)$ if and only if $\Ext_R^i(R/\fm,R)=0$ for all $i>d$. We now extend this characterization to relative Gorenstein modules
over arbitrary rings. To do so, we first introduce the following definition:

\begin{definition}\label{3.1} Let $\fa$ be a proper ideal of $R$ and $M$ an $R$-module. The $\fa$-{\it relative injective dimension}
of $M$ is defined by $\fa-\id_RM=\sup \{i\in \mathbb{N}_0 \mid \Ext_R^i(R/\fa,M)\neq 0\}.$
\end{definition}

It is worth mentioning that $\fm-\id_RM=\id_RM$ for any finitely generated module $M$ over a local ring $(R,\fm)$.

\begin{theorem}\label{3.2} Let $\fa$ be an ideal of $R$, $M$ a finitely generated $R$-module with $M\neq \fa M$ and set $c=\cd(\fa,R)$.
Then the following are equivalent:
\begin{enumerate}
\item[(i)] $M$ is $\fa$-relative Gorenstein;
\item[(ii)] $M$ is $\fa$-relative maximal Cohen-Macaulay and $\fa-\id_R\left({\RH}_{\fa}^c\left(M\right)\right)=0$;
\item[(iii)] $M$ is $\fa$-relative maximal Cohen-Macaulay and $\Ext_R^i(R/\fa,M)=0$ for all $i>c$.
\end{enumerate}
\end{theorem}

\begin{prf} Let $N$ be an $\fa$-relative maximal Cohen-Macaulay $R$-module. Then ${\RH}_{\fa}^i\left(N\right)=0$ for all $i\neq c$.
Thus by \cite[Proposition 2.1]{RZ2}, one has $$\Ext_R^i(R/\fa,\RH_{\fa}^c(N))\cong \Ext_R^{i+c}(R/\fa,N) \  \  \   (*)$$ for all
$i\geq 0$.

(i) $\Rightarrow$ (ii) As $M\neq \fa M$ and $M$ is $\fa$-relative Gorenstein, by Theorem \ref{2.19}(iii), we see that $M$ is
$\fa$-relative maximal Cohen-Macaulay, and $\Ext_R^c(R/\fa,M)$ is the only non-vanishing Ext module of $R/\fa$ and $M$. Thus
from $(*)$, we conclude that $\fa-\id_R\left({\RH}_{\fa}^c\left(M\right)\right)=0$.

(ii) $\Rightarrow$ (i) Since $M$ is $\fa$-relative maximal Cohen-Macaulay, we have $\grade(\fa,M)=c$. Hence, $\Ext_R^i(R/\fa,M)=0$
for all $i<c$. On the other hand, as $\fa-\id_R\left({\RH}_{\fa}^c\left(M\right)\right)=0$, by $(*)$, we conclude that
$\Ext_R^i(R/\fa,M)=0$ for all $i>c$. Therefore, $M$ is $\fa$-relative Gorenstein.

(i) $\Leftrightarrow$ (iii) is obvious by Theorem \ref{2.19}(iii) and the definition.
\end{prf}

Recall that over a local ring $R$, a Gorenstein module is a maximal Cohen-Macaulay module of finite injective dimension.
The fact that for a finitely generated $R$-module $M$ and a maximal ideal $\fm$ of $R$, all local cohomology modules
${\RH}_{\fm}^i\left(M\right)$ are Artinian, enables us to record the following corollary.

\begin{corollary}\label{3.3} Let $(R,\fm)$ be a local ring of dimension $d$ and $M$ a nonzero finitely generated $R$-module.
Then the following are equivalent:
\begin{enumerate}
\item[(i)] $M$ is  Gorenstein.
\item[(ii)] $M$ is maximal Cohen-Macaulay and ${\RH}_{\fm}^d\left(M\right)$ is an Artinian injective $R$-module.
\end{enumerate}
\end{corollary}

In order to present the next result, we need the following definition; see \cite[Definition 3.1]{PDR}.

\begin{definition}\label{3.4} Let $\fa$ be a proper ideal of $R$ and $c=\cd(\fa,R)$. Assume that $R$ is $\fa$-relative
Cohen-Macaulay. The $\fa$-{\it relative dualizing module} of $R$ is defined by $$\Omega_{\fa}=\Hom_R(\RH_{\fa}^c(R),
\underset{\fm\in \Max R}\bigoplus \text{E}_R(R/\fm)).$$
\end{definition}

This section's final main result is now ready to be presented.

\begin{theorem}\label{3.5} Let $R$ be a complete semi-local ring, $\fa$ an ideal of $R$ contained in its Jacobson radical and
set $c=\cd(\fa,R)$. Assume that $R$ is $\fa$-relative Cohen-Macaulay and $\Omega_{\fa}$ can be identified with an ideal $\fb$ of
$R$. If $\RH_{\fa}^c(R/\fb)=0$, then the ring $R/\fb$ is $(\fa+\fb)/\fb$-relative Gorenstein.
\end{theorem}

\begin{prf} Set $T=R/\fb$. From the exact sequence $$0\lo \fb\lo R\lo T\lo 0, \  \  \  \ (*)$$ we deduce the exact sequence
$$\cdots \lo \RH_{\fa}^i(R)\lo \RH_{\fa}^i(T)\lo \RH_{\fa}^{i+1}(\fb)\lo \cdots  .$$ By \cite[Theorem 3.5(ii)]{PDR},  $\RH_{\fa}^{i}(\fb)\cong\RH_{\fa}^{i}(\Omega_{\fa})=0$ for all $i\neq c$. So, as $R$ is $\fa$-relative Cohen-Macaulay
and $\RH_{\fa}^c(T)=0$, from the above exact sequence, we conclude that $\RH_{\fa}^i(T)=0$ for all $i\neq c-1$. As $T\neq \fa T$,
it turns out that $T$ is $\fa T$-relative Cohen-Macaulay and $\cd(\fa T,T)=c-1$.  Hence, $\Omega_{\fa T}\cong \Ext_R^{1}(T,\fb)$
by \cite[Theorem 3.7]{PDR}.

Next, by \cite[Theorem 3.5(iii)]{PDR}, we have $$\Hom_R(\fb,\fb)\cong \underset{\fm\in \Max R}\prod \widehat{R_{\fm}}
\cong R.$$ In particular, $\Ann_R(\fb)=0$. This implies that $\Hom_R(T,\fb)=0$. So, applying the functor $\Hom_R(-,\fb)$
to the exact sequence $(*)$ yields the exact sequence $$0\lo \fb \lo R\lo \Ext_R^{1}(T,\fb)\lo 0.$$ Thus, $$T\cong
\Ext_R^{1}(T,\fb)\cong \Omega_{\fa T}.$$ Therefore, \cite[Lemma 3.3(ii)]{PDR} implies that $$\Ext_T^i(T/\fa T,T)\cong
\Ext_T^i(T/\fa T,\Omega_{\fa T})=0$$ for all $i\neq c-1$, and so the ring $T$ is $\fa T$-relative Gorenstein.
\end{prf}

According to Bass's theorem, a local ring $(R,\fm)$ is Cohen-Macaulay if it possesses a nonzero finitely generated module
of finite injective dimension. Moreover, a local ring $(R,\fm)$ is Gorenstein if it possesses a nonzero cyclic $R$-module
of finite injective dimension. This might lead us to guess that a ring $R$ is $\fa$-relative Cohen-Macaulay (resp.
$\fa$-relative Gorenstein) if it admits a nonzero finitely generated (resp. cyclic) module $M$ of finite $\fa$-relative
injective dimension. As we will see in the following example, this is not the case.

\begin{example}\label{3.6} Let $\Bbbk$ be a field, $S=\Bbbk[x_1,x_2,y_1,y_2]$ and $\fm=\langle x_1, x_2, y_1, y_2
\rangle$. Let $\fa$ be the edge ideal of the cycle graph $C_4$, given in Example \ref{2.20}. So, $\fa=\langle
x_1x_2,x_2y_1,y_1y_2,y_2x_1 \rangle$. We observed in Example \ref{2.20} that $\fa=\langle x_1, y_1 \rangle \cap
\langle x_2, y_2 \rangle$. This yields the exact sequence $$0\lo S/\fa\lo S/\langle x_1, y_1 \rangle
\oplus S/\langle x_2, y_2 \rangle \lo S/\fm \lo 0, \ \  (*)$$
which immediately implies that $\dim_S(S/\fa)=2$. Taking into account the long exact sequence of local cohomology modules
induced by $(*),$ yields that $\RH_{\fm}^0(S/\fa)=0$ and $\RH_{\fm}^1(S/\fa)\cong S/\fm\neq 0$. Hence, $\depth_S(S/\fa)=1$.
Thus, by \cite[Theorem 1]{Ly} and the Auslander-Buchsbaum formula, we deduce that $$\cd(\fa,S)=\pd_S(S/\fa)=3.$$ As $S$ is
Cohen-Macaulay, we get $$\grade(\fa,S)=\dim S-\dim_S(S/\fa)=2.$$ Thus, $S$ is not $\fa$-relative Cohen-Macaulay. On the
other hand, we have $\fa-\id_SS\leq \pd_S(S/\fa)=3$.
\end{example}

We need the following lemma for our next example.

\begin{lemma}\label{3.7a} Let $(R,\fm)$ be a local ring and $\fa$ a  proper ideal of $R$ such that $\pd_R(R/\fa)<\infty$.
Then $\fa-\id_RM=\pd_R(R/\fa)$ for every nonzero finitely generated $R$-module $M$.
\end{lemma}

\begin{prf} Let $M$ be a nonzero finitely generated $R$-module. As $\pd_R(R/\fa)<\infty$, by \cite[\S 7, Lemma 1(iii)]{M}, it turns
out that $$\pd_R(R/\fa)=\sup \{i\in \mathbb{N}_0 \mid \Ext_R^{i}(R/\fa,N)\neq 0\}$$ for every nonzero finitely generated $R$-module
$N$. In particular, $\pd_R(R/\fa)=\fa-\id_RM$.
\end{prf}

When $(R,\fm)$ is a $d$-dimensional local ring, vanishing of the Ext modules $\Ext_R^i(R/\fm,R)$ for all $i>d$ implies that $R$ is
Gorenstein. This might suggest that if $\fa$ is a proper ideal of $R$ and $M$ a finitely generated $R$-module such that
$\Ext_R^i(R/\fa,M)=0$ for all $i>\cd(\fa,R)$, then $M$ is $\fa$-relative Gorenstein. While this is not the case by Example \ref{3.6},
we also give the following simpler example.

\begin{example}\label{3.8b} Let $(R,\fm)$ be a local ring and $x_1\in \fm$ a nonzero-divisor on $R$. Set $\fa=\langle x_1 \rangle$.
Then $\pd_R(R/\fa)=1$. Let $M$ be any nonzero finitely generated $R$-module which is annihilated by $x_1$. Then, by the choice of
$M$ and Lemma \ref{3.7a}, we have $\Ext_R^i(R/\fa,M)\neq 0$ for $i=0,1$. So, although $\Ext_R^i(R/\fa,M)=0$ for all $i>1=\cd(\fa,R)$,
the $R$-module $M$ is not $\fa$-relative Gorenstein.
\end{example}

Our next example requires the following lemma.

\begin{lemma}\label{3.9c} Let $\fa$ an ideal of $R$ and $M$ a finitely generated $R$-module with $M\neq \fa M$. Let $x\in \fa$ be
a nonzero-divisor on $M$. Then $\fa-\id_RM=\fa-\id_R(M/xM)$.
\end{lemma}

\begin{prf} Set $\overline{M}=M/xM$. The short exact sequence $$0\lo M\overset{x}\lo M\lo \overline{M}\lo 0$$ yields the
exact sequence $$\cdots \lo \Ext_R^{i-1}(R/\fa,\overline{M})\lo \Ext_R^{i}(R/\fa,M)\overset{x}\lo
\Ext_R^{i}(R/\fa,M)\lo \Ext_R^{i}(R/\fa,\overline{M})\lo $$ $$\lo \Ext_R^{i+1}(R/\fa,M)\overset{x}\lo \Ext_R^{i+1}(R/\fa,M)\lo\cdots .
\      \  (\dag)$$ Let $n$ be a non-negative integer. Assume that $\Ext_R^{j}(R/\fa,M)=0$ for all $j>n$. Then from $(\dag)$, we conclude
that $\Ext_R^{j}(R/\fa,\overline{M})=0$ for all $j>n$. Hence, $\fa-\id_R\overline{M}\leq \fa-\id_RM$.

Next, assume that $\Ext_R^{j}(R/\fa,\overline{M})=0$ for all $j>n$. Then for every $j>n$, from $(\dag)$, we see that the zero
map $$\Ext_R^{j}(R/\fa,M)\overset{x}\lo \Ext_R^{j}(R/\fa,M)$$ is surjective, and so $\Ext_R^{j}(R/\fa,M)=0$. Thus, $\fa-\id_RM\leq \fa-\id_R\overline{M}$. Therefore, $\fa-\id_RM=\fa-\id_R(M/xM)$.
\end{prf}

Let $\fa$ be an ideal of $R$ and $M$ an $\fa$-relative Gorenstein $R$-module. One may guess that if $x\in \fa$ is a nonzero-divisor on
$M$, then the $R$-module $M/xM$ is also $\fa$-relative Gorenstein. The following example demonstrates that this is not true.

\begin{example}\label{3.10d} Let $c\geq 1$ be an integer. Let $\fa$ be an ideal of $R$ that is generated by an $R$-regular sequence
$x_1, x_2, \dots, x_c$. Then, clearly,  $R$ is $\fa$-relative Gorenstein. Now, Lemma \ref{3.9c} implies that $$\fa-\id_R(R/\langle x_1
\rangle)=\fa-\id_RR,$$ while $$\grade(\fa,R/\langle x_1 \rangle)=\grade(\fa,R)-1=\fa-\id_RR-1.$$ Hence, the $R$-module $R/\langle x_1
\rangle$ is not $\fa$-relative Gorenstein.
\end{example}

A local ring $(R,\fm)$ is Gorenstein if and only if $R$ is Cohen-Macaulay and it possesses an irreducible parameter ideal.
Accordingly, for a proper ideal $\fa$ of $R$, we may expect that $R$ is $\fa$-relative Gorenstein if and only if $R$ is
$\fa$-relative Cohen-Macaulay and there exists an $\fa$-s.o.p $x_1, x_2,\dots, x_c$ such that the ideal $\langle x_1, x_2,
\dots, x_c\rangle$ is irreducible. However, as the following example shows, this is not the case.

\begin{example}\label{3.11e} Let $(R,\fm)$ be a $d$-dimensional Cohen-Macaulay local ring which is not Gorenstein. Since $R$ is not
Gorenstein, it has no irreducible parameter ideal. Let $x_1, x_2, \dots, x_d\in \fm$ be a system of parameters of $R$. As $R$ is
Cohen-Macaulay, $x_1, x_2, \dots, x_d$ is an $R$-regular sequence. Set $\fa=\langle x_1, x_2, \dots, x_d\rangle$. Then $R$ is
$\fa$-relative Gorenstein, while there is no irreducible ideal generated by an $\fa$-s.o.p of $R$.
\end{example}

For a local ring $(R,\fm)$ and a nonzero finitely generated $R$-module $M$ of finite injective dimension, it is known that
$\id_RM=\depth R$. As a consequence, we may conjecture that if $\fa$ is an ideal of $R$ and $M$ a finitely generated $R$-module
with $M\neq \fa M$ such that $\fa-\id_RM<\infty$, then $\fa-\id_RM=\grade(\fa,R)$. The following example shows that this is not
the case as well.

\begin{example}\label{3.12} Let $(R,\fm)$ be a regular local ring and $\fa$ a nonzero proper ideal of $R$ such that $\cd(\fa,R)\neq
\pd_R(R/\fa)$. Then Lemma \ref{3.7a} implies that $\fa-\id_RR\neq \cd(\fa,R)$. More precisely, let $\Bbbk$ be a field and consider
the formal power series ring $R=\Bbbk[[x,y]]$. Then, $R$ is a regular local ring with the unique maximal ideal $\fm=\langle x,y
\rangle$. Let $\fa=\langle xy,x^2 \rangle$. One has $$\Ass_R(R/\fa)=\{\langle x \rangle,\fm\}.$$
As $\fm\in \Ass_R(R/\fa)$, it follows that $\depth_R(R/\fa)=0$, and so by Lemma \ref{3.7a} and the Auslander-Buchsbaum formula, we get
that $$\fa-\id_RR=\pd_R(R/\fa)=2.$$ On the other hand, $$\Rad(\fa)=\langle x \rangle \cap \fm=\langle x
\rangle.$$ Thus, we have
$$\begin{array}{ll}
1&\leq \grade(\fa,R)\\
&\leq \cd(\fa,R)\\
&=\cd(\langle x \rangle,R)\\
&=1.
\end{array}$$
Therefore, $R$ is $\fa$-relative Cohen-Macaulay and $\grade(\fa,R)=1$. Consequently, $\fa-\id_RR\neq \grade(\fa,R)$.
\end{example}

\begin{definition}\label{3.13}  A finitely generated $R$-module $C$ is called {\it semidualizing} if it satisfies the following conditions:
\begin{itemize}
\item[(i)] the homothety map $\chi_C^R:R \lo \Hom_R\left(C,C\right)$ is an isomorphism, and
\item[(ii)] $\Ext^i_R\left(C,C\right)=0$ for all $i>0$.
\end{itemize}
\end{definition}

\begin{definition}\label{3.14} Let $C$ be a semidualizing module of $R$.
\begin{itemize}
\item[(i)] The {\it Auslander class} $\mathscr{A}_C\left(R\right)$ is the class of all $R$-modules $M$ for which the natural
map	$\gamma_M^C:M\lo \Hom_R\left(C,C\otimes_RM\right)$ is an isomorphism, and $$\Tor^R_i\left(C,M\right)=0=\Ext_R^i\left(C,C\otimes_RM\right)$$ for all	$i\geq 1$.
\item[(ii)] The {\it Bass class} $\mathscr{B}_C\left(R\right)$ is the class of all $R$-modules $M$ for which the evaluation map
$\xi_M^C:C\otimes_R\Hom_R\left(C,M\right)\lo M$ is an isomorphism, and $$\Ext^i_R\left(C,M\right)=0=\Tor^R_i\left(C,\Hom_R\left(C,M\right)\right)$$ for all
$i\geq 1$.
\end{itemize}
\end{definition}

\begin{notation}\label{3.15} Let $\fa$ be a proper ideal of $R$ and $n$ a non-negative integer.
\begin{itemize}
\item[(i)] $\text{CM}^n_{\fa}(R)$ stands for the full subcategory of $\fa$-relative Cohen-Macaulay $R$-modules $M$
with $\cd\left(\fa,M\right)=n$.
\item[(ii)] $\mathscr{G}_{\mathfrak{a}}\left(R\right)$ stands for the full subcategory of $\fa$-relative Gorenstein $R$-modules.
\end{itemize}
\end{notation}

Let $\fa$ be a proper ideal of $R$, $C$ a semidualizing module of $R$ and $n$ a non-negative integer. By \cite[Theorem 6.3]{PDR},
there is an equivalence of categories:
\begin{displaymath}
\xymatrix{\mathscr{A}_C\left(R\right)\bigcap \text{CM}^n_{\fa}(R) \ar@<0.7ex>[rrr]^-{C\otimes_R-} &
{} & {} & \mathscr{B}_C\left(R\right)\bigcap \text{CM}^n_{\fa}(R).  \ar@<0.7ex>[lll]^-{\Hom_R\left(C,-\right)}}
\end{displaymath}

The following example illustrates that $\text{CM}^n_{\fa}(R)$ cannot be replaced by $\mathscr{G}_{\fa}(R)$ in the above equivalence
of categories.

\begin{example}\label{3.16} Let $(R,\fm)$ be a non-Gorenstein Cohen-Macaulay local ring with a dualizing module $\omega_R$. Then
$\omega_R$ is a semidualizing module of $R$ and it is $\fm$-relative Gorenstein. On the other hand, it is easy to see that $\omega_R\in
\mathscr{B}_{\omega_R}\left(R\right)$. Hence $\omega_R\in \mathscr{B}_{\omega_R}\left(R\right)\bigcap \mathscr{G}_{\fm}(R)$, while
$$\Hom_R(\omega_R,\omega_R)\cong R\notin \mathscr{A}_{\omega_R}\left(R\right)\bigcap \mathscr{G}_{\fm}(R).$$ Thus, the functors
\begin{displaymath}
\xymatrix{\mathscr{A}_{\omega_R}\left(R\right)\bigcap \mathscr{G}_{\fm}(R) \ar@<0.7ex>[rrr]^-{\omega_R\otimes_R-} &
{} & {} & \mathscr{B}_{\omega_R}\left(R\right)\bigcap \mathscr{G}_{\fm}(R)  \ar@<0.7ex>[lll]^-{\Hom_R\left(\omega_R,-\right)}}
\end{displaymath}
do not induce an equivalence of categories.
\end{example}

Let $(R,\fm)$ be a Cohen-Macaulay local ring with a dualizing module $\omega_R$ and $\text{MCM}(R)$ denote the full subcategory
of maximal Cohen-Macaulay $R$-modules. There is a well-known duality of categories:
\begin{displaymath}
\xymatrix{\text{MCM}(R) \ar@<0.7ex>[rrr]^-{\Hom_R(-,\omega_R)} &
{} & {} &  \text{MCM}(R).  \ar@<0.7ex>[lll]^-{\Hom_R(-,\omega_R)}}
\end{displaymath}

One might expect that $\text{MCM}(R)$ can be replaced by $\mathscr{G}_{\fa}(R)$ in the above duality of categories. However, the
following example shows that this is not true as well.

\begin{example}\label{3.17} Let $(R,\fm)$ be a non-Gorenstein Cohen-Macaulay local ring with a dualizing module $\omega_R$. Then
$\omega_R$ is $\fm$-relative Gorenstein, while $\Hom_R(\omega_R,\omega_R)\cong R$ is not $\fm$-relative Gorenstein. Hence, the
functor
\begin{displaymath}
\xymatrix{\mathscr{G}_{\fm}(R) \ar@<0.7ex>[rrr]^-{\Hom_R(-,\omega_R)} &
{} & {} &  \mathscr{G}_{\fm}(R)  \ar@<0.7ex>[lll]^-{\Hom_R(-,\omega_R)}}
\end{displaymath}
does not induce a duality of categories.
\end{example}

\section{Relative regular modules}

Theorem \ref{4.1a} generalizes Lemma \ref{2.18} to relative regular modules, provided that $\fa$ is contained in the Jacobson radical
of $R$. A characterization of relative regular modules is given in Theorem \ref{4.4d}. We also establish the fact that every relative
regular module is relative Gorenstein in Proposition \ref{4.6f}. In Propositions \ref{4.8h} and \ref{4.9i}, we demonstrate that relative
regular modules remain invariant under certain equivalences and dualities of categories.

\begin{theorem}\label{4.1a} Let $\fa$ be an ideal of $R$ contained in its Jacobson radical and $M$ a nonzero $\fa$-relative regular
$R$-module. Then $\pd_R(M/\fa M)=\pd_RM+\cd\left(\fa,R\right)$.
\end{theorem}

\begin{prf} Set $c=\cd\left(\fa,R\right)$. As $M$ is nonzero and $\fa$ is contained in the Jacobson radical of $R$, it follows that
$M\neq \fa M$. So, $\fa$ has generators $x_1, x_2, \dots, x_c$ which form both an $M$-regular sequence and an $R$-regular sequence.
In particular, $\fa=0$ if and only if $c=0$. Since the claim holds trivially if $\fa=0$, we may assume that $c>0$. By induction
on $c$, we may assume that $c=1$. Set $\overline{R}=R/\langle x_1 \rangle$ and $\overline{M}=M/\langle x_1 \rangle M$. Then
$\overline{M}$ is an $\fa \overline{R}$-relative regular $\overline{R}$-module. For every $R$-module $N$, from the short exact
sequence $$0\lo M\overset{x_1}\lo M\lo \overline{M}\lo 0,$$  we get the following exact sequence $$\cdots \lo \Ext_R^{i}(M,N)
\overset{x_1}\lo \Ext_R^{i}(M,N)\lo \Ext_R^{i+1}(\overline{M},N)\lo \Ext_R^{i+1}(M,N)\lo \cdots.  \  \   \  (*)$$
Let $n$ be a non-negative integer. If $\pd_RM=n$, then from $(*)$, we deduce that $\Ext_R^{j}(\overline{M},N)=0$ for all $j>n+1$.
Thus, $$\pd_R(\overline{M})\leq \pd_RM+1.$$ Next,  suppose that $\pd_R(\overline{M})=n$. Then for every $j>n-1$, from $(*)$, we
conclude that the map $$\Ext_R^{j}(M,N)\overset{x_1}\lo \Ext_R^{j}(M,N)$$ is surjective, and so $\Ext_R^{j}(M,N)=0$ by Nakayama's
lemma. Hence, $$\pd_RM\leq \pd_R(\overline{M})-1.$$ Therefore, $$\pd_R(\overline{M})=\pd_RM+1,$$ as required.
\end{prf}

To present the next main result, we need the following two lemmas.

\begin{lemma}\label{4.2b} Let $\fa$ be an ideal of $R$ contained in its Jacobson radical. If $R$ is $\fa$-relative regular, then every $\fa$-relative maximal Cohen-Macaulay $R$-module is $\fa$-relative regular.
\end{lemma}

\begin{prf} Set $c=\cd(\fa,R)$. Assume that $R$ is $\fa$-relative regular. So, there is an $R$-regular sequence $x_1, \dots, x_c$
that generates $\fa$. Let $M$ be an $\fa$-relative maximal Cohen-Macaulay $R$-module.  Then, by Theorem \ref{2.19}(iv), it turns out
that $\cd(\fa,M)=c$. Now, from the definition, it is obvious that $x_1, \dots, x_c$ is an $\fa$-s.o.p of $M$.
Therefore, $x_1,\dots, x_c$ is an $M$-regular sequence by Lemma \ref{2.7b}, and so $M$ is $\fa$-relative regular.
\end{prf}

In what follows, $\mu(\fa)$ stands for the minimum number of generators of an ideal $\fa$.

\begin{lemma}\label{4.3c} Let $\fa$ be a proper ideal of $R$. Then $R$ is $\fa$-relative regular if and only if $\grade(\fa,R)=\mu(\fa)$.
\end{lemma}

\begin{prf} Set $c=\cd(\fa,R)$. First, assume that $R$ is $\fa$-relative regular. Then $\fa$ can be generated by an $R$-regular sequence
$x_1,\dots, x_c$. So, $c\leq \grade(\fa,R)\leq c$ and $\mu(\fa)\leq c$. On the other hand, as $\RH_{\fa}^{i}\left(M\right)=0$ for all
$i>\mu(\fa)$, it follows that $c\leq \mu(\fa)$. Thus, $\grade(\fa,R)=\mu(\fa)$.

Conversely, suppose that $\grade(\fa,R)=\mu(\fa)$. Then $$\mu(\fa)=\grade(\fa,R)\leq c\leq \mu(\fa),$$ and so $\grade(\fa,R)=c$. Now, \cite[Exercise 1.2.21]{BH} yields that $\fa$ can be generated by an $R$-regular sequence of length $c$. Thus, $R$ is $\fa$-relative
regular.
\end{prf}

\begin{theorem}\label{4.4d} Let $\fa$ be an ideal of $R$ contained in its Jacobson radical and $M$ a nonzero finitely generated $R$-module.
Then the following are equivalent:
\begin{enumerate}
\item[(i)] $M$ is $\fa$-relative regular;
\item[(ii)]  $\grade(\fa,M)=\grade(\fa,R)=\mu(\fa)$.
\end{enumerate}
\end{theorem}

\begin{prf} Set $c=\cd(\fa,R)$.

(i) $\Rightarrow$ (ii) Assume that $M$ is $\fa$-relative regular. Then $\fa$ can be generated by an $M$-regular sequence
$x_1, \dots, x_c$ of length $c$. Hence, $$c\leq \grade(\fa,M)\leq \cd(\fa,M)\leq \mu(\fa)\leq c,$$ and so $\grade(\fa,M)
=\mu(\fa)$. Since $M$ is $\fa$-relative regular, by the definition, it follows that $R$ is also $\fa$-relative regular.
So, $\grade(\fa,R)=\mu(\fa)$ by Lemma \ref{4.3c}.

(ii) $\Rightarrow$ (i) Assume that $$\grade(\fa,M)=\grade(\fa,R)=\mu(\fa).$$ The equality $\grade(\fa,R)=\mu(\fa)$ yields
that $R$ is $\fa$-relative regular by Lemma \ref{4.3c}. Hence, $R$ is $\fa$-relative Cohen-Macaulay by Theorem \ref{2.19},
and so $$\grade(\fa,M)=\grade(\fa,R)=c.$$ Thus $M$ is $\fa$-relative maximal Cohen-Macaulay, and so $M$ is $\fa$-relative
regular by Lemma \ref{4.2b}.
\end{prf}

A local ring $(R,\fm)$ is regular if and only if $\dim R=\pd_R(R/\fm)$. So, one may guess that if $\fa$ is a proper ideal of
$R$, then $R$ is $\fa$-relative regular if and only if $\cd\left(\fa,R\right)=\pd_R(R/\fa)$. The following example indicates
that this is not the case.

\begin{example}\label{4.5e}  Let $\Bbbk$ be a field, $S=\Bbbk[x_1,x_2,y_1,y_2]$ and $\fa=\langle x_1x_2,x_2y_1,y_1y_2,y_2x_1
\rangle$. In Example \ref{3.6}, we observed that $$\cd(\fa,S)=\pd_S(S/\fa)=3.$$ On the other hand, $\fa$ can not be generated
by an $S$-regular sequence of length 3, because $\grade(\fa,S)=2$. Hence, $S$ is not $\fa$-relative regular.
\end{example}

Next, we show that every $\fa$-relative regular $R$-module is $\fa$-relative Gorenstein.

\begin{proposition}\label{4.6f} Let $\fa$ be a proper ideal of $R$ and $M$ an $\fa$-relative regular $R$-module. Then
$$\Ext_R^i(R/\fa,M)\cong
\begin{cases} M/\fa M \  \   \hspace{.5cm} \text{if} \hspace{.5cm}  i=\cd(\fa,R) \\
0  \hspace{1.7cm}  \text{if} \hspace{.5cm}  i\neq \cd(\fa,R).
\end{cases}
$$
In particular, $M$ is $\fa$-relative Gorenstein.
\end{proposition}

\begin{prf} Set $c=\cd(\fa,R)$. Clearly, we may assume that $M\neq \fa M$. Hence, the ideal $\fa$ has generators $x_{1}, x_2,
\dots, x_{c}$, which form both an $M$-regular sequence and an $R$-regular sequence. As $x_{1}, x_2, \dots, x_{c}$ is an
$M$-regular sequence, we see $$c\leq \grade(\fa,M)\leq \cd(\fa,M)\leq c,$$ and so $\grade(\fa,M)=c$. In particular,
$\Ext_R^i(R/\fa,M)=0$ for all $i<c$. From the definition, it follows that $R$ is $\fa$-relative regular, and so Lemma
\ref{2.18} implies that $\pd_R(R/\fa)=c$. So, $\Ext_R^i(R/\fa,M)=0$ for all $i>c$. Thus $M$ is $\fa$-relative Gorenstein.

Finally, \cite[Lemma 1.2.4]{BH} yields that
$$
\begin{array}{ll}
\Ext_R^c(R/\fa,M)&\cong \Hom_R(R/\fa,M/\langle x_{1}, x_2, \dots, x_{c}\rangle M)\\
&=\Hom_R(R/\fa,M/\fa M)\\
&\cong M/\fa M.\\
\end{array}
$$
\end{prf}

\begin{notation}\label{4.7g} Let $\fa$ be a proper ideal of $R$. Let $\mathscr{R}_{\mathfrak{a}}\left(R\right)$ denote the full
subcategory of $\fa$-relative regular $R$-modules.
\end{notation}

\begin{proposition}\label{4.8h} Let $\fa$ be a proper ideal of $R$. Assume that $R$ is $\fa$-relative regular. For every
semidualizing module $C$ of $R$, there is an equivalence of categories:
\begin{displaymath}
\xymatrix{\mathscr{A}_C\left(R\right)\cap \mathscr{R}_{\mathfrak{a}}\left(R\right) \ar@<0.7ex>[rrr]^-{C\otimes_R-} &
{} & {} & \mathscr{B}_C\left(R\right)\cap \mathscr{R}_{\mathfrak{a}}\left(R\right).  \ar@<0.7ex>[lll]^-{\Hom_R\left(C,-\right)}}
\end{displaymath}
\end{proposition}

\begin{prf} Set $c=\cd(\fa,R)$. Because of the equivalence
\begin{displaymath}
\xymatrix{\mathscr{A}_C\left(R\right) \ar@<0.7ex>[rrr]^-{C\otimes_R-} &
{} & {} & \mathscr{B}_C\left(R\right),  \ar@<0.7ex>[lll]^-{\Hom_R\left(C,-\right)}}
\end{displaymath}
it is enough to show that a finitely generated $R$-module $N\in \mathscr{A}_C\left(R\right)$ belongs to
$\mathscr{R}_{\mathfrak{a}}\left(R\right)$ if and only if $C\otimes_RN$  belongs to $\mathscr{R}_{\mathfrak{a}}\left(R\right)$.

Let $M\in \mathscr{A}_C\left(R\right)$ be a finitely generated $R$-module. As $R$ is $\fa$-relative regular, the ideal $\fa$ can
be generated by an $R$-regular sequence $\underline{x}=x_1, \dots, x_c$. By \cite[Lemma 3.2(ii)]{AbDT}, $\underline{x}$
is an $M$-regular sequence if and only if it is a $C\otimes_RM$-regular sequence. Thus, $M\in \mathscr{R}_{\mathfrak{a}}
\left(R\right)$ if and only if $C\otimes_RM\in \mathscr{R}_{\mathfrak{a}}\left(R\right)$.
\end{prf}

\begin{proposition}\label{4.9i}  Let $\mathfrak{J}$ denote the Jacobson radical of a complete semi-local ring $R$. Assume that $R$
is $\mathfrak{J}$-relative regular. Then there is a duality of categories:
\begin{displaymath}
\xymatrix{\mathscr{R}_{\mathfrak{J}}\left(R\right) \ar@<0.7ex>[rrr]^-{\Hom_R\left(-,\Omega_{\mathfrak{J}}\right)} &
{} & {} & \mathscr{R}_{\mathfrak{J}}\left(R\right).  \ar@<0.7ex>[lll]^-{\Hom_R\left(-,\Omega_{\mathfrak{J}}\right)}}
\end{displaymath}
\end{proposition}

\begin{prf} Let $\text{MCM}_{\mathfrak{J}}(R)$ stands for the full subcategory of $\mathfrak{J}$-relative maximal Cohen-Macaulay
$R$-modules. Proposition \ref{4.6f} and Theorem \ref{2.19}(iii) yield that $\mathscr{R}_{\mathfrak{J}}\left(R\right)\setminus \{0\}
\subseteq \text{MCM}_{\mathfrak{J}}(R)$. On the other hand, as $R$ is $\mathfrak{J}$-relative regular, Lemma \ref{4.2b} implies that
$\text{MCM}_{\mathfrak{J}}(R)\subseteq \mathscr{R}_{\mathfrak{J}}\left(R\right)\setminus \{0\}$. Thus $\mathscr{R}_{\mathfrak{J}}
\left(R\right)=\text{MCM}_{\mathfrak{J}}(R)\cup \{0\}$, and so the claim follows by \cite[Corollary 5.3]{PDR}.
\end{prf}


\end{document}